\documentclass[12pt]{amsart}
\usepackage{a4wide,color}
\usepackage{amsmath,amssymb}
\usepackage[active]{srcltx}

\allowdisplaybreaks

\let\pa\partial

\let\eps\varepsilon

\newcommand{\R}{{\mathbb R}}
\newcommand{\T}{{\mathbb T}}

\newtheorem{theorem}{Theorem}

\newtheorem{proposition}[theorem]{Proposition}

\let\ga=\gamma

\let\eps=\varepsilon

\let\la=\lambda

\let\pa=\partial

\let\ka=\kappa

\let\Om=\Omega

\begin{document}

\title[Gross-Pitaevskii Equation]{Anelastic Approximation of the Gross-Pitaevskii equation for
General Initial Data}

\author[C.-K. Lin]{Chi-Kun Lin}
\address{Chi-Kun Lin, Department of Mathematical Sciences, Xian Jiaotong Liverpool University, Suzhou, China}
\email{cklin@math.nctu.edu.tw}

\author[K.-C. Wu]{Kung-Chien Wu}
\address{Kung-Chien Wu, Department of Mathematics,
National Cheng Kung University, Tainan, Taiwan }
\email{kcwu@mail.ncku.edu.tw; kungchienwu@gmail.com}

\date{\today}

\begin{abstract}
We perform a rigorous analysis of  the anelastic approximation for the Gross-Pitaevskii equation with $x$-dependent chemical potential.
For general initial data and periodic boundary condition,
we show that as $\eps\to 0$, equivalently the Planck constant tends to zero, the density $|\psi^{\eps}|^{2}$ converges toward the chemical potential $\rho_{0}(x)$ and the velocity field converges
to the anelastic system. When the chemical potential is a constant, the anelastic system will reduce to the incompressible Euler equations.
The resonant effects the singular limit process and it can be overcome because of oscillation-cancelation.
\end{abstract}



\maketitle

\section{Introduction}
We consider a superfluid governed by the Gross-Pitaevskii equation
\cite{[PS]}
\begin{equation}\label{int.1.a11}
i\hbar\pa_{t}\psi+\frac{\hbar^{2}}{2m}\Delta\psi -(V_0|\psi|^{2}-
E)\psi=0\,,
\end{equation}
appropriate to a weakly interacting Bose gas. This is a nonlinear
Schr\"odinger equation for a single-particle wave function
$\psi(x,t)$ governing an assembly of bosons of mass $m$, with $V_0$
the strength of the $\delta$-function interaction potential between
the bosons, $E$ the chemical potential, and $\hbar$ the Planck
constant. In this paper, we will assume the unit strength $V_0=1$
and the chemical potential $E$ is a function of position, i.e.,
$E=\rho_0(x)$;
\begin{equation}\label{int.1.a1}
i\hbar\pa_{t}\psi+\frac{\hbar^{2}}{2m}\Delta\psi
-\big(|\psi|^{2}-\rho_{0}(x)\big)\psi=0\,.
\end{equation}
In order to investigate the singular limit, we introduce the scaled variables $\widetilde{t}=\eps t$ and $\widetilde{x}=x$ with $0<\eps\ll 1$, and the Planck constant is also rescaled as $\hbar =\eps^{1+\alpha}$, where $0<\alpha <\infty$.
Then after dropping the tilde, (\ref{int.1.a1}) becomes
\begin{equation}\label{int.1.a}
i\eps^{\alpha}\pa_{t}\psi^{\eps}+\frac{\eps^{2\alpha}}{2}\Delta\psi^{\eps}
-\frac{1}{\eps^{2}}\big(|\psi^{\eps}|^{2}-\rho_{0}(x)\big)\psi^{\eps}=0\,.
\end{equation}
The initial condition is complemented by
\begin{equation}
\psi^{\eps}(x,0)=\psi^{\eps}_0(x)\,,\quad x\in \Omega\,.
\end{equation}

The Gross-Pitaevskii equation (\ref{int.1.a}) defines a Hamiltonian
system. The conserved Hamiltonian is a Ginzburg-Landau energy,
namely
\begin{equation}
\frac{1}{2}\int \eps^{2\alpha} |\nabla\psi^{\eps}|^2
+\frac{1}{\eps^{2}}\big( |\psi^{\eps}|^2-\rho_{0}(x)\big)^{2}dx=C_1\,.
\end{equation}
Similarly, the momentum or current
\begin{equation}
\frac{i}{2}\eps^{\alpha}\int
\Big(\psi^{\eps}\nabla\overline{\psi^{\eps}}
-\overline{\psi^{\eps}}\nabla\psi^{\eps}\Big)dx=C_2
\end{equation}
is formally conserved. Here $\overline{\psi^{\eps}}$ denotes the complex conjugate of the wave function $\psi^\eps$. Another quantity which is formally conserved
by the flow is the mass or charge
\begin{equation}
\int |\psi^{\eps}|^2 dx =\int
\psi^{\eps}\overline{\psi^{\eps}} dx=C_3\,.
\end{equation}
According to the above conservative quantities we have the
hydrodynamical variables: density (or charge) $\rho^{\eps}$,
momentum (or current) $J^{\eps}$, and energy $e^{\eps}$ given
respectively by
\begin{equation}\begin{array}{c}\label{int.1.c}
\displaystyle\rho^{\eps}=|\psi^{\eps}|^2,\qquad
J^{\eps}=\frac{i}{2}\eps^{\alpha}
\big(\psi^{\eps}\nabla\overline{\psi^{\eps}}
-\overline{\psi^{\eps}}\nabla\psi^{\eps}\big)\,,
\\ \\
\displaystyle
e^{\eps}=\frac{1}{2}\eps^{2\alpha} |\nabla\psi^{\eps}|^2
+\frac{1}{2\eps^{2}}\big( |\psi^{\eps}|^2-\rho_{0}\big)^{2}.
\end{array}\end{equation}
Under the above definition, we can define the initial conditions as follows: initial charge $\rho_{0}^{\eps}(x)=\rho^\eps(x,0)=|\psi^{\eps}_0(x)|^2$ and
initial momentum (current) $J_{0}^{\eps}(x)=J_0^\eps(x,0)$. The
associated local conservations of mass(charge), momentum(current)
and energy in terms of hydrodynamic variables are given by
\begin{equation}\label{int.2.a1}
\partial_t\rho^{\eps}+ \nabla\cdot J^{\eps}=0\,,
\end{equation}
\begin{equation}\begin{array}{lll}\label{int.2.b1}
\displaystyle\partial_t J^{\eps}+\nabla\cdot
\bigg(\frac{J^{\eps}\otimes J^{\eps}}{\rho^{\eps}} \bigg) +\frac{1}{\eps^{2}}\rho^{\eps}\nabla (\rho^{\eps}-\rho_{0})
=\frac{\eps^{2\alpha}}{4}\nabla\cdot
\Big(\rho^{\eps}\nabla^2\log\rho^{\eps}\Big)\,,
\end{array}\end{equation}
\begin{equation}\begin{array}{l}\label{int.2.c1}
\displaystyle\partial_t e^{\eps}
+\nabla\cdot\bigg(e^{\eps}\frac{J^{\eps}}{\rho^{\eps}}
+\frac{1}{\eps^{2}}J^{\eps}(\rho^{\eps}-\rho_{0})\bigg)
=\frac{\eps^{2\alpha}}{4}\nabla\cdot
\bigg(\frac{J^{\eps}}{\rho^{\eps}}\Delta \rho^{\eps} -\nabla\cdot
J^{\eps}\frac{ \nabla \rho^{\eps}}{\rho^{\eps}}\bigg)\,.
\end{array}
\end{equation}
Equations (\ref{int.2.a1})--(\ref{int.2.b1}) comprise a closed system governing
$\rho^\eps$ and $J^\eps$, which have the form of a perturbation
of the compressible Euler equations with an extra potential
$\rho_0(x)$. Formally, letting $\eps\to 0$, the uniform boundedness of energy
$e^{\eps}$ and strict convexity of the potential will imply $\rho^{\eps}\to \rho_{0}$. We also assume $J^\eps \to \rho_0 v$ for some proper function $v$. The limit equations will then be the anelastic system
\begin{equation}\begin{array}{c}\label{hs.4.cc}
\displaystyle\partial_t(\rho_{0} v) +\hbox{div} \big(\rho_{0} v\otimes
v\big) +\rho_{0}\nabla \pi =0\,,
\\[2mm]
v(x,0)=v_0(x)\,,\quad \nabla\cdot (\rho_{0} v)=0\,,
\end{array}\end{equation}
where the pressure $\pi$ is the formal limit of
$\frac{1}{\eps^{2}}(\rho^{\eps}-\rho_{0})$. When $\rho_0$ is a constant, the anelastic system (\ref{hs.4.cc}) will reduce to the incompressible Euler equation,
\begin{equation}\begin{array}{c}\label{hs.4.ccc}
\partial_t v + v\cdot \nabla v\ +\nabla \pi =0\,,
\\[2mm]
v(x,0)=v_0(x)\,,\quad \nabla\cdot v=0\,.
\end{array}
\end{equation}
The reader is referred to \cite{[Lions96]} for the detail study of the incompressible Euler equations.
For $x\in \Omega\subset \R^2$, (\ref{hs.4.cc}) can be served as the lake equations which may be seen as the low Froude number limit of the usual inviscid shallow water equations when the initial height converges to a nonconstant function depending on the space variable. Here we have to interpret $\pi(x,t)$ as the surface height variation (see \cite{[BGL], [Bresch-Metivier], [LMT]} and references therein).
For the viscous shallow water equations and the convergence to the quasi-geostrophic model we will refer to \cite{[BD2003]}.
Note that the chemical potential $\rho_0(x)$ plays the role of depth of the basin at location $x$ in the lake equations.
However, the above formal discussion of the singular limit $\eps \to 0$ can not be made directly in
(\ref{int.2.a1})--(\ref{int.2.b1}) since $\log \rho^\eps$ may be undefined.
Indeed, we can write the dispersive term on the right hand side of
(\ref{int.2.b1}) in different ways:
\begin{equation}\label{int.1.aaa}
\frac{1}{4}\nabla\cdot\left(\rho^\eps\nabla^2\log\rho^\eps\right)
=\frac{1}{2} \rho^\eps\nabla
\left(\frac{\Delta\sqrt{\rho^\eps}}{\sqrt{\rho^\eps}}\right)=
\frac{1}{4} \Delta \nabla \rho^\eps - \nabla\cdot \left(\nabla
\sqrt{\rho^\eps}\otimes \nabla \sqrt{\rho^\eps}\right)\,.
\end{equation}
Instead of the momentum and energy equations
(\ref{int.2.b1})--(\ref{int.2.c1}), we prefer to represent them as
\begin{align}\label{int.2.b}
 \partial_t J^{\eps}
 &+\displaystyle\frac{1}{2}\eps^{2\alpha} \nabla\cdot
 \Big((\nabla\psi^{\eps}\otimes\nabla\overline{\psi^\eps}
 +\nabla\overline{\psi^\eps}\otimes\nabla\psi^{\eps})
 -\nabla^{2}(|\psi^{\eps}|^{2})\Big)\nonumber\\
&+\frac{1}{2\eps^{2}}\nabla(\rho^{\eps}-\rho_{0})^{2}
+\frac{1}{\eps^{2}}\rho_{0}\nabla(\rho^{\eps}-\rho_{0})=0\,,
\end{align}
and
\begin{equation}\begin{array}{ccc}\label{int.2.c}
\displaystyle\partial_t e^{\eps}
-\frac{1}{2}\eps^{2\alpha}\nabla\cdot\Big(
 \nabla\psi^{\eps}\partial_{t}\overline{\psi^{\eps}}
+\nabla\overline{\psi^{\eps}}\partial_t\psi^{\eps}\Big)=0\,,
\end{array}
\end{equation}
which are more suitable for the estimates of the modulated energy as we shall see in the next section.

Besides the Wigner transform (see \cite{[GLM]} for the comparison with the WKB analysis), there are two other ways to study the semiclassical limit of the Schr\"odinger equation. In the framework of
Sobolev space and a defocusing nonlinearity, E. Grenier \cite{[Gr]}
introduce the modified Madelung transformation to rewrite the
nonlinear Schr\"odinger equation as a {\it linear} perturbation of
quasilinear symmetric hyperbolic system.

Another approach to prove
the convergence of nonlinear Schr\"odinger equation to Euler
equations came from Y. Brenier \cite{[Bre]}, who proved the convergence of the Vlasov-Poisson system to the incompressible Euler equations, following an ideal due to P.-L. Lions in \cite{[Lions96]} (see also \cite{[LinWu08S]} for Klein-Gordon equation and \cite{[LZ]} for Gross-Pitaevskii equation).
The main advantage of
this approach is that we compare the two solutions  in more or less
the energy space for nonlinear Schr\"odinger equation.
It was immediately extended by Masmoudi in \cite{[Mas01]} to general initial conditions allowing the presence of high oscillation in time
and this idea is also applied to quantum hydrodynamical model of semiconductor \cite{[LiLin]} and the Schr\"odinger-Poisson system in Coulomb gauge \cite{[LinWongWu11]}. Recently, D. Han-Kwan combines the modulated energy and relative entropy method to prove the quasineutral limit of the Vlasov-Poisson system with massless electrons \cite{[Han-Kwan]}. In this paper we will apply the modulated energy method to prove the convergence from the Gross-Pitaevskii equation (\ref{int.1.a11})
to the anelastic system (\ref{hs.4.cc}).

The question of anelastic-type limits in fluid dynamics has received considerable attention recently. As is well-known, an anelastic approximation
is a filtering approximation for the equations of motion that
eliminates sound waves by assuming that the flow has velocities and
phase speeds much smaller than the speed of sound. This
approximation has been used to model astrophysical and geophysical
fluids \cite{[Maj03], [Ped]}. The rigorous derivation of the anelastic limit starting from the diffusive systems have been recently studied from a mathematical point of view in \cite{[BGL], [FMNS-2008],  [Masmoudi07], [Wu12]}. We will refer to \cite{[Bresch-Metivier-2010]} for the compressible Euler-type systems and \cite{[Gallagher]} for the study of the asymptotic behavior of a fluid submitted to a strong external magnetic field.

As has been pointed out in the study of the low Mach number limit of fluid dynamics, the main problem is to investigate the average effect of fast acoustic waves on the slow incompressible motion \cite{[BDGL], [LionsM], [MS], [Schochet94]}. For the anelastic limit, it is the fast acoustic wave for inhomogeneous media. The heterogeneity of the medium is modeled through the chemical potential $\rho_0(x)$. But, this fast wave equation is independent of the solution and using this property, the anelastic limit can be carried out rigorously.
The main difficulty in this limit is the crossing eigenvalues phenomena, thus, we have to study the resonant effect of the oscillating part $\mathcal{Q}_{2}(V,V)$ (see Proposition \ref{pr3}).
The reader is referred to \cite{[Mas01a]} for constant coefficient acoustic wave equation and \cite{[Bresch-Metivier-2010], [MS]} for non-constant case.

The rest of the paper is organized as follows. In section 2, in addition to the main theorem, we also present the weighted Helmholtz decomposition which will play the essential role for the anelastic limit. It was introduced by P-L Lions in \cite{[Lions96]} to study the density-dependent Navier-Stokes equations. Later on Bresch et al. \cite{[BDGL]} applied this decomposition to prove the low Mach number limit of the polytropic viscous flows. Recently, Feireisl et al. \cite{[FMNS-2008]} successfully applied the weighted Helmholtz decomposition, combining the spectral theorem to derive the anelastic approximation of the compressible Navier-Stokes equations.

Section 3 is devoted to the proof of the main theorem. Besides the wave group introduced by Schochet \cite{[Schochet94]} (see also \cite{[BDGL], [LiLin],  [LionsM], [Mas01]}), we also need the spectral analysis of the associated highly oscillating wave operator \cite{[FMNS-2008]}. We employ the modulated energy method to prove the convergence. Since the limiting density is not a constant, we have to introduce two correction terms of the modulated energy which measure the variation of the density. Indeed, we will show they tend to 0 as $\eps\to 0$.
\section{Main Theorem and Function Space}

Since $\rho_0$ is strictly positive, we can define the weighted
space as follows. Let $\Omega\subset \R^n$ be open and $\sigma(x)={1\over \rho_0(x)}$ be the weighted function, then the weighted
square integrable space $L^2({dx\over
\rho_0})=L_\sigma^{2}(\Omega)$ consists of all measurable
functions $f$ that satisfy
$$
\int_\Omega |f(x)|^2 \sigma(x) dx <\infty\,.
$$
The resulting $L_\sigma^{2}(\Omega)$-norm of $f$ is defined by
\begin{equation}
\|f\|_{L_\sigma^{2}(\Omega)}= \left(\int_\Omega |f(x)|^2 \sigma(x) dx
\right)^{1/2}.
\end{equation}
The space $L_\sigma^{2}(\Omega)$ is naturally equipped with
the following inner product:
\begin{equation}
\langle f,g\rangle_\sigma\equiv\int_\Omega f\cdot
g\sigma(x)dx\,.
\end{equation}
The weighted Sobolev space $H^1_\sigma(\Omega)$ consists of all functions $f$ with weak derivatives $Df$ satisfying
\begin{equation}
\|f\|_{H^1_\sigma(\Omega)}= \left(\int_\Omega \big(|f(x)|^2 + |Df(x)|^2\big)\sigma(x)dx \right)^{1/2}<\infty.
\end{equation}
To derive the anelastic system we need to introduce the space $\mathbb{H}_{\rho_{0}}[L^2_\sigma(\Omega)]$ defined by
\begin{equation}
\mathbb{H}_{\rho_{0}}[L^2_\sigma(\Omega)]
=\{ \rho_{0}v\in L^2_\sigma(\Omega)\,|\, \hbox{div } (\rho_{0}v)=0\}\,.
\end{equation}
Since $\mathbb{H}_{\rho_{0}}[L^2_\sigma]$ is a closed subspace of $L^2_\sigma$, then by projection theorem the Hilbert space $L^2_\sigma(\Omega)$ admits an orthogonal projection
\begin{equation}\label{2.5}
L^2_\sigma(\Omega) =\mathbb{H}_{\rho_{0}}[L^2_\sigma(\Omega)]\oplus
\mathbb{H}^{\perp}_{\rho_{0}}[L^2_\sigma(\Omega)]\,,
\end{equation}
where the orthogonal complement $\mathbb{H}^{\perp}_{\rho_{0}}[L^2_\sigma(\Omega)]$ is given by
\begin{equation}
\mathbb{H}^\perp_{\rho_{0}}[L^2_\sigma(\Omega)]
=\left\{ \rho_0\nabla \Psi\,|\, \Psi \in H^1_\sigma(\Omega), \, \int_\Omega \Psi dx =0\right\}\,.
\end{equation}
In the sequel, we will consider the periodic domain $\Omega=\T^n$.
Associated with the orthogonal projection (\ref{2.5}) is the
weighted Helmholtz decomposition in the form \cite{[FMNS-2008]}
\begin{equation}
\mathbb{H}_{\rho_{0}}f= f- \rho_{0}\nabla\Psi\,,\qquad
\mathbb{H}^{\perp}_{\rho_{0}}f=\rho_{0}\nabla\Psi\,,
\end{equation}
where $\Psi\in H^1_\sigma(\T^{n})$ is the unique solution of the
problem
\begin{equation}
\int_{\T^{n}}\rho_{0}\nabla\Psi\cdot\nabla\eta dx
=\int_{\T^{n}}f\cdot\nabla\eta dx\,,\qquad \forall \eta\in
H^1_\sigma(\T^{n})\,.
\end{equation}
That is, $\Psi$ is a weak solution of the uniformly elliptic partial differential equation satisfying the zero mean condition
\begin{equation}
\nabla\cdot(\rho_{0}\nabla\Psi)=\nabla\cdot f\,, \quad \int_{\T^{n}}\Psi dx=0\,.
\end{equation}
Note that $\mathbb{H}_{\rho_{0}}$ is an orthogonal projection on $L^2_\sigma$, and the two projectors $\mathbb{H}_{\rho_{0}}f$ and $\mathbb{H}^{\perp}_{\rho_{0}}f$ are orthogonal with respect to the inner product
$\langle \cdot,\cdot\rangle_{\sigma}$.
The reader is referred to \cite{[BDGL], [FMNS-2008], [Lions96]} for
mathematical properties of such projectors.
In order to make sure that (\ref{int.1.a}) and the limit systems (\ref{pe.0.a}) and (\ref{pe.3.a}) are well-defined, we shall from now on
impose the following conditions on $\psi^{\eps}_{0}$, $\rho_0$, $J_0$, $v_0$ and $w_0$:
\bigskip

\noindent{\rm(A1)}\quad $\psi^{\eps}_{0}\in H^{{n\over2}+3}(\mathbb{T}^{n};\mathbb{C})$, this will guarantee the local existence and uniqueness of classical solution of the Gross-Pitaevskii equation (\ref{int.1.a}).

\bigskip

\noindent{\rm(A2)}\quad  $\rho_{0}\in C^{s}(\mathbb{T}^{n})$, $s>\frac{n}{2}+1$, $\rho_{0}\geq c>0$, the initial kinetic, potential and quantum energies satisfy
\begin{equation}\begin{array}{l}\label{int.2.d}
\displaystyle \frac{1}{\sqrt{\rho^{\eps}_{0}}}J^{\eps}_{0}
\to\frac{1}{\sqrt{\rho_{0}}}J_{0} \quad \hbox{in}\quad
L^{2}{(\T^{n})}\,,
\\ \\
\displaystyle\frac{\rho_{0}^{\eps}-\rho_{0}}{\eps}
\to \varphi_{0} \quad \hbox{in}\quad L^{2}{(\T^{n})}\,,
\\ \\
\eps^{\alpha}\nabla\sqrt{\rho_{0}^{\eps}}\to 0  \quad \hbox{in}\quad
L^{2}{(\T^{n})}\,.
\end{array}\end{equation}

\bigskip

\noindent{\rm(A3)}\quad
$J_{0}=\mathbb{H}_{\rho_{0}}J_{0}\oplus\mathbb{H}^{\perp}_{\rho_{0}}J_{0}
=\rho_{0}v_{0}+\rho_{0}\nabla w_{0}$, where $\sqrt{\rho_{0}}v_{0}\in
H^{s}(\T^{n})$ and  $( \varphi_{0},\sqrt{\rho_{0}}\nabla w_{0})\in
H^{s}(\T^{n})$ for $s>\frac{n}{2}+3$. This condition will guarantee the
local existence and uniqueness of smooth solution of the
anelastic system (\ref{pe.0.a}) and the oscillating part
(\ref{pe.3.a}) (see Proposition \ref{pro77}).
\bigskip

The main result of this paper is stated as follows.
\begin{theorem}\label{maintheo2}
Let $\alpha>0$ and  $\psi^{\eps}$ be the solution of the Gross-Pitaevskii equation
$(\ref{int.1.a})$ with $\psi^\eps_{0}$ satisfying the assumptions ${\rm (A1)}$--${\rm(A3)}$. There then exists
$T>0$ such that
\begin{equation}\label{ZDL.4.a}
\rho^{\eps}\to \rho_{0} \quad\hbox{strongly in}\quad
L^{\infty}\big([0,T];{L^{2}(\T^{n})}\big)\,,
\end{equation}
\begin{equation}\label{ZDL.4.b}
J^{\eps}\rightharpoonup \rho_{0}v \quad \hbox{weakly $*$ in} \quad
L^{\infty}\big([0,T];L^{4/3}(\T^{n})\big)\,,
\end{equation}
where $v$ solves the anelastic system
\begin{equation}\label{pe.0.a}
\left\{\begin{array}{l} \pa_{t}(\rho_{0}v)+\nabla\cdot (\rho_{0} v\otimes v)+\rho_{0}\nabla\pi=0\,,
\\[2mm]
v(x,0)=v_{0}(x)\,,\quad \nabla\cdot (\rho_{0}v)=0\,.
\\
\end{array}
\right.\end{equation}
In particular, if $\rho_0=1$ or any constant, then $v$ will be the solution of the incompressible Euler equations
\begin{equation}\label{pe.0.b}
\left\{\begin{array}{l} \pa_{t}v+v\cdot \nabla v+\nabla\pi=0\,,
\\[2mm]
v(x,0)=v_{0}(x)\,,\quad \nabla\cdot v=0\,.
\\
\end{array}
\right.\end{equation}

\end{theorem}

\section{Proof of the Main Theorem}

We divide the proof into several steps:

\noindent{\bf Step 1.} {\it Spectral analysis of the wave group.\it} Analogous to the low Mach number limit in fluid dynamics we consider the perturbation of $\rho^\eps$ near the equilibrium $\rho_0$
\begin{equation}\label{pe.1.a}
\rho^{\eps}=\rho_{0}+\eps\varphi^{\eps}\,,
\end{equation}
i.e., $\varphi^\eps$ is the density fluctuation.
By weighted Helmholtz decomposition the current $J^\eps$ can be rewritten as
$$
J^{\eps}=\mathbb{H}_{\rho_{0}}J^{\eps}\oplus
         \mathbb{H}^{\perp}_{\rho_{0}}J^{\eps}
        =\rho_{0}u^{\eps}+\rho_{0}\nabla w^{\eps}\,.
$$
Applying the operator $\mathbb{H}^{\perp}_{\rho_{0}}$ to the momentum equation (\ref{int.2.b}) and using the weighted incompressibility $\nabla \cdot (\rho_{0}u^{\eps})=0$, the equations (\ref{int.2.a1}) and (\ref{int.2.b}) can be rewritten in terms of $\varphi^\eps$ and $w^\eps$ as
\begin{equation}\label{pe.1.b}
\left\{\begin{array}{l} \displaystyle\pa_{t}\varphi^{\eps}+\frac{1}{\eps}\hbox{div}(\rho_{0}\nabla w^{\eps})=0\,,
\\[2mm]
\displaystyle\pa_{t}(\sqrt{\rho_{0}}\nabla w^{\eps})+\frac{1}{\eps}\sqrt{\rho_{0}}\nabla\varphi^{\eps}=\frac{1}{\sqrt{\rho_{0}}} F^{\eps}\,,
\\
\end{array}
\right.\end{equation}
where
\begin{align}\label{pe.1.c}
F^{\eps}=&-\frac{\eps^{2\alpha}}{2} \mathbb{H}^{\perp}_{\rho_{0}}\nabla\cdot
 \big(\nabla\psi^{\eps}\otimes\nabla\overline{\psi^{\eps}}
 +\nabla\overline{\psi^{\eps}}\otimes\nabla\psi^{\eps}
\big)\nonumber\\
&-\frac{1}{2}\mathbb{H}^{\perp}_{\rho_{0}}\nabla(\varphi^{\eps})^{2}
+\frac{\eps^{2\alpha}}{4}\mathbb{H}^{\perp}_{\rho_{0}}\nabla\Delta\rho^{\eps}\,.
\end{align}
It is obvious from (\ref{pe.1.b}) that $\pa_{t}\varphi^{\eps}$ and $\pa_{t}(\sqrt{\rho_{0}}\nabla w^{\eps})$ are of order $O(1/\eps)$ and are highly oscillatory as $\eps\to 0$. Therefore we have to introduce the wave group in order to filter out the fast oscillating wave.
Let $\mathcal{L}(\tau)=e^{\tau L}$, $\tau\in \R$ be the evolution group associated with the operator $L$ which is defined, according to (\ref{pe.1.b}), on $\mathbb{D}=L^{2}(\T^{n})\times\{h=\sqrt{\rho_{0}}\nabla w: h \in L^{2}(\T^{n})\}$ by
\begin{equation}\label{pe.1.d}
L\left(
\begin{array}{c}
\phi \\ \sqrt{\rho_{0}}\nabla w
\end{array}
\right)=-\left(
\begin{array}{c}
\hbox{div} (\rho_{0}\nabla w) \\ \sqrt{\rho_{0}}\nabla\phi
\end{array}
\right)
\,.
\end{equation}
The eigenvalue problem associated with $L$ is
\begin{equation}\label{pe.1.d}
L\left(
\begin{array}{c}
\phi \\ \sqrt{\rho_{0}}\nabla w
\end{array}
\right)=
-\left(
\begin{array}{c}
\hbox{div} (\rho_{0}\nabla w) \\ \sqrt{\rho_{0}}\nabla\phi
\end{array}
\right)
=\la\left(
\begin{array}{c}
\phi \\ \sqrt{\rho_{0}}\nabla w
\end{array}
\right)
\,.
\end{equation}
The second row gives $\la \nabla w=-\nabla\phi$, hence the eigenvalue problem of $L$ is equivalent to the eigenvalue problem of the uniformly elliptic operator
\begin{equation}\label{pe.1.d.1}
-\nabla\cdot(\rho_{0}\nabla\phi)=-\la^{2}\phi\,.
\end{equation}
For convenience we define the operator $\mathcal{A}$ by $\mathcal{A}\phi\equiv -\nabla\cdot(\rho_{0}\nabla\phi)$.
Since $\rho_0$ is positive and bounded away from zero, the associated bilinear form of $-\nabla\cdot(\rho_{0}\nabla\phi)$ is coercive, $\mathcal{A}$ is a linear unbounded operator and its inverse operator $\mathcal{A}^{-1}$ is self-adjoint and compact by Rellich lemma.
Thus there exists an orthonormal basis $\{\chi_j\}_{j=1}^\infty$ of $L^2(\T^n)$, i.e.,
$$
\int_{\T^n} \chi_l(x)\cdot \chi_j(x) dx =\delta_{lj}
$$
and a sequence $\{\ka_j\}_{j=1}^\infty$ such that
$$
0<\ka_1\le \ka_2\le \ka_3 \le \cdots  ,\qquad \ka_j\to \infty
$$
$$
\mathcal{A}\chi_j = \ka_j \chi_j\,,\qquad \forall j \,,
$$
where each eigenvalue is repeated according to its multiplicity (which is known to be finite).
If we pull back to the original operator $L$,
then we have $\la_{j}^{\pm}=\pm i\sqrt{\ka_{j}}$. Thus the eigenvalues and the associated eigenfunctions of $L$ are given by
\begin{equation}
\bigg\{i\sqrt{\ka_{j}}\,,\left(
\begin{array}{c}
\chi_{j} \\ \frac{i}{\sqrt{\ka_{j}}}\sqrt{\rho_{0}}\nabla\chi_{j}
\end{array}
\right)\bigg\}
\quad\hbox{and}\quad
\bigg\{-i\sqrt{\ka_{j}}\,,\left(
\begin{array}{c}
\chi_{j} \\ \frac{-i}{\sqrt{\ka_{j}}}\sqrt{\rho_{0}}\nabla\chi_{j}
\end{array}
\right)\bigg\}\,.
\end{equation}
Since $V\in \mathbb{D}$ is real, we can represent it as the eigenfunction expansion
\begin{equation}\label{3.8}
V=\sum_{j=1}^{\infty}a_{j}^{+}\left(
\begin{array}{c}
\chi_{j} \\ \frac{i}{\sqrt{\ka_{j}}}\sqrt{\rho_{0}}\nabla\chi_{j}
\end{array}
\right)+a_{j}^{-}\left(
\begin{array}{c}
\chi_{j} \\ \frac{-i}{\sqrt{\ka_{j}}}\sqrt{\rho_{0}}\nabla\chi_{j}
\end{array}
\right)
\end{equation}
where $a_{j}^{+}$ is the complex conjugate of $a_{j}^{-}$, i.e., $a_{j}^{+}=\overline{a_{j}^{-}}$.  Accordingly, by spectral decomposition of $L$, we have
\begin{equation}\label{3.9}
\begin{array}{l}
\displaystyle\mathcal{L}\Big(\frac{t}{\eps}\Big)V=\displaystyle\sum_{j=1}^{\infty}a_{j}^{+}e^{i\sqrt{\ka_{j}}t/\eps}\left(
\begin{array}{c}
\chi_{j} \\ \frac{i}{\sqrt{\ka_{j}}}\sqrt{\rho_{0}}\nabla\chi_{j}
\end{array}
\right)+a_{j}^{-}e^{-i\sqrt{\ka_{j}}t/\eps}\left(
\begin{array}{c}
\chi_{j} \\ \frac{-i}{\sqrt{\ka_{j}}}\sqrt{\rho_{0}}\nabla\chi_{j}
\end{array}
\right)\,.
\end{array}\end{equation}
For convenience, we denote $\mathcal{L}_{1}(\tau)$ the first component and $\mathcal{L}_{2}(\tau)$ the rest $n$ components of $\mathcal{L}(\tau)$, hence
\begin{equation}\label{3.10}
\begin{array}{l}
\eps\pa_{t}\mathcal{L}_{1}(\frac{t}{\eps})V=-\hbox{div}(\sqrt{\rho_{0}}\mathcal{L}_{2}(\frac{t}{\eps})V)\,,
\\[2mm]
\eps\pa_{t}\mathcal{L}_{2}(\frac{t}{\eps})V=-\sqrt{\rho_{0}}\nabla(\mathcal{L}_{1}(\frac{t}{\eps})V)\,.
\end{array}\end{equation}

\begin{proposition}
The operator $\mathcal{L}(\tau)$ is isometry on
$L^{2}(\T^{n})$ and is bounded on $H^{s}(\T^{n})$ for all $s\in\R
$ and for all $\tau\in \R$.
\end{proposition}
\noindent{\it Proof.} For any $V\in \mathbb{D}$, if $V\in L^{2}(\T^{n})$, (note that the eigenvalues of $L$ are pure imaginary numbers) we have
$$
\int_{\T^{n}}|\mathcal{L}(\tau)V|^{2}dx=\int_{\T^{n}}|V|^{2}dx\,.
$$
For $V\in H^{s}(\T^{n}), s>0$ and any $|\alpha|\leq s$,
$$
\int_{\T^{n}}|D^{\alpha}\mathcal{L}(\tau)V|^{2}dx\leq\|V\|^{2}_{H^{s}(\T^{n})}\,,
$$
and if $V\in H^{-s}(\T^{n}), s>0$, for any $f\in \mathbb{D}\cap
C_{0}^{\infty}(\T^{n})$
\begin{align*}
\int_{\T^{n}}\mathcal{L}(\tau)V(x)\cdot f(x) dx&=\int_{\T^{n}}V(x)\cdot \mathcal{L}(-\tau)f(x) dx\\
&  \leq\|\mathcal{L}(-\tau)f\|_{H^{s}(\T^{n})}\leq
\|f\|_{H^{s}(\T^{n})}\,,
\end{align*}
so $\mathcal{L}(\tau)$ is unitary in $L^{2}(\T^{n})$ and is bounded in
Hilbert space $H^{s}(\T^{n})$ for all $s\in\R$ and $\tau>0$. \qed

In the sequel, we shall denote
$$
U^{\eps} = \left(
\begin{array}{c}
\varphi^{\eps} \\  \sqrt{\rho_{0}}\nabla w^{\eps}
\end{array}
\right)\,,\qquad
V^{\eps} = \mathcal{L}\Big(\frac{-t}{\eps}\Big)\left(
\begin{array}{c}
\varphi^{\eps} \\  \sqrt{\rho_{0}}\nabla w^{\eps}
\end{array}
\right)\,.
$$
With this notation, we can rewrite (\ref{pe.1.b}) as
$$
\pa_t U^{\eps}={1\over \eps} LU^{\eps} +\frac{1}{\sqrt{\rho_{0}}}\widehat{F}^{\eps}\,,
$$
or equivalently
\begin{equation}\label{pe.1.g}
\pa_{t}V^{\eps}=\mathcal{L}\Big(\frac{-t}{\eps}\Big)\frac{1}{\sqrt{\rho_{0}}}\widehat{F}^{\eps}\,,
\end{equation}
where $\widehat{v}$ denotes the column vector $(0,v)^{t}$.
\begin{proposition} Let $m\geq \frac{n}{2}+5$.
There exists a subsequence of $\{V^\eps\}_\eps$ which we still
denote by $\{V^\eps\}_\eps$ and $\overline{V}\in
L^{1}\big([0,T];H^{-m}(\T^{n})\big)$, such that
\begin{equation}\label{ZDL.6.f}
V^{\eps}\to \overline{V} \quad\hbox{strongly in} \quad
L^{1}\big([0,T];H^{-m}(\T^{n})\big)\,.
\end{equation}
\end{proposition}
\noindent{\it Proof.} For $m\geq \frac{n}{2}+5$, it is easy to see
that
$$
\hbox{$U^{\eps}$ is uniformly bounded in
$L^{\infty}\big([0,T],L^{2}(\T^{n})\big)$}
$$
and
$$
\hbox{$\frac{1}{\sqrt{\rho_{0}}}F^{\eps}$ is uniformly bounded in
$L^{2}\big([0,T];H^{-m}(\T^{n})\big)$}\,.
$$
By the boundedness of ${\mathcal{L}}$, we can show that
$$
\hbox{$V^{\eps}$ is uniformly bounded in
$L^{\infty}\big([0,T],L^{2}(\T^{n})\big)$}
$$
and
$$
\hbox{$\pa_{t}V^{\eps}$ is uniformly bounded in
$L^{2}\big([0,T];H^{-m}(\T^{n})\big)$.}
$$
Therefore, we deduce from the Lions-Aubin's lemma that there exists
a subsequence of $\{V^{\eps}\}_{\eps}$ which we still denote by
$\{V^{\eps}\}_{\eps}$ and $\overline{V}\in
L^{2}\big([0,T];H^{-m}(\T^{n})\big)$ such that
\begin{equation}\label{pe.1.h}
V^{\eps}\to \overline{V} \quad\hbox{strongly in} \quad
L^{1}\big([0,T];H^{-m}(\T^{n})\big)\,.
\end{equation}
\qed

\noindent{\bf Step 2.} {\it Analysis of the oscillating equation.\it}
Let $V, V_{1}, V_{2}\in \mathbb{D}\cap
L^{2}(\T^{n})$, and $u$ be any vector in $\R^{n}$ such that
$\hbox{div}(\rho_{0}u)=0$, we define the linear form $B_1$ by
\begin{equation}\begin{array}{l}\label{3.14}
B_{1}(u,V)=\hbox{div}\Big(\sqrt{\rho_{0}}u\otimes
\mathcal{L}_{2}(\frac{t}{\eps})V+\sqrt{\rho_{0}}
\mathcal{L}_{2}(\frac{t}{\eps})V\otimes u\Big)\,
\end{array}\end{equation}
and the bilinear form $B_2$ by
\begin{align}\label{3.15}
B_{2}(V_{1},V_{2})&=\frac{1}{2}\hbox{div}\Big(\mathcal{L}_{2}\Big(\frac{t}{\eps}\Big)V_{1}\otimes
\mathcal{L}_{2}\Big(\frac{t}{\eps}\Big)V_{2}+\mathcal{L}_{2}\Big(\frac{t}{\eps}\Big)V_{2}
\otimes \mathcal{L}_{2}\Big(\frac{t}{\eps}\Big)V_{1}\Big)\\
&\phantom{xx}{}+\frac{1}{2}\nabla\Big(\mathcal{L}_{1}\Big(\frac{t}{\eps}\Big)V_{1}\cdot
\mathcal{L}_{1}\Big(\frac{t}{\eps}\Big)V_{2}\Big)\,.\nonumber
\end{align}
We can rewrite (\ref{pe.1.g}) as
\begin{align}\label{pe.1.g11}
\displaystyle\pa_{t}V^{\eps}&=\mathcal{L}\Big(\frac{-t}{\eps}\Big)\left(
\begin{array}{c}
0 \\  \frac{1}{\sqrt{\rho_{0}}}\mathbb{H}^{\perp}_{\rho_{0}}B_{1}(u^{\eps},V^{\eps})
+\frac{1}{\sqrt{\rho_{0}}}\mathbb{H}^{\perp}_{\rho_{0}}B_{2}(V^{\eps},V^{\eps})
\end{array}
\right)\\
&\phantom{xx}{}+\frac{\eps^{2\alpha}}{2} \mathcal{L}\Big(\frac{-t}{\eps}\Big)\left(
\begin{array}{c}
0 \\ \frac{1}{\sqrt{\rho_{0}}}\mathbb{H}^{\perp}_{\rho_{0}}\rho^\eps\nabla
\left(\frac{\Delta\sqrt{\rho^\eps}}{\sqrt{\rho^\eps}}\right)
\end{array}
\right)\,.\nonumber
\end{align}
If we had sufficient compactness in space, then we could pass to the limit as $\eps\to 0$ in (\ref{pe.1.g11}) and obtain the following system for the oscillating parts:
\begin{equation}\label{pe.1.i}
\pa_{t}\overline{V}+\mathcal{Q}_{1}(v,\overline{V})+\mathcal{Q}_{2}(\overline{V},\overline{V})=0\,,
\end{equation}
where $v$ is the strong solution of the anelastic system (\ref{pe.0.a}).
Here $\mathcal{Q}_{1}$ and $\mathcal{Q}_{2}$ are respectively a linear and a bilinear forms of $V$ defined by
\begin{equation}\label{pe.1.j}
\mathcal{Q}_{1}(u,V)=\lim_{\tau\to \infty}\frac{1}{\tau}\int_{0}^{\tau}\mathcal{L}(-s)
\left(
\begin{array}{c}
0 \\  \frac{1}{\sqrt{\rho_{0}}}\mathbb{H}^{\perp}_{\rho_{0}}B_{1}(u,V)
\end{array}
\right)ds
\end{equation}
and
\begin{equation}\label{pe.1.ja}
\mathcal{Q}_{2}(V,V)=\lim_{\tau\to \infty}\frac{1}{\tau}\int_{0}^{\tau}\mathcal{L}(-s)
\left(
\begin{array}{c}
0 \\  \frac{1}{\sqrt{\rho_{0}}}\mathbb{H}^{\perp}_{\rho_{0}}B_{2}(V,V)
\end{array}
\right)ds
\end{equation}
for any divergence free vector fields $u\in L^{2}(\T^{n})$ and any
$V=(\phi, \sqrt{\rho_{0}}\nabla q)^{t}\in L^{2}(\T^{n})$.
Actually we have the following proposition:

\begin{proposition}\label{pr1}
Let $1/r_{1}+1/r_{2}=1$.
For all $ u\in L^{r_{1}}([0,T];L^{2}(\T^{n}))$
and $V\in L^{r_{2}}([0,T];L^{2}(\T^{n}))$, we have the following
weak converges
\begin{equation}\begin{array}{c}\label{pe.2.a}
{\displaystyle w\hbox{-}\lim_{\eps\to 0}}
\mathcal{L}(\frac{-t}{\eps}) \left(
\begin{array}{c}
0 \\  \frac{1}{\sqrt{\rho_{0}}}\mathbb{H}^{\perp}_{\rho_{0}}B_{1}(u,V)
\end{array}
\right)=\mathcal{Q}_{1}(u,V)
\end{array}\end{equation}
and
\begin{equation}\begin{array}{c}\label{pe.2.a}
{\displaystyle w\hbox{-}\lim_{\eps\to 0}}
\mathcal{L}(\frac{-t}{\eps}) \left(
\begin{array}{c}
0 \\  \frac{1}{\sqrt{\rho_{0}}}\mathbb{H}^{\perp}_{\rho_{0}}B_{2}(V,V)
\end{array}
\right)=\mathcal{Q}_{2}(V,V)\,.
\end{array}\end{equation}
\end{proposition}

  One can show the above converges by using the almost-periodic
functions and the reader is referred to \cite{[Mas01a]} and the references therein. Using the symmetry of $\mathcal{Q}_{2}$, we also have the following proposition.
\begin{proposition}\label{pr4} Let $s\in \R$ and $1/p+1/q=1$. The following identity
\begin{equation}\begin{array}{c}
{\displaystyle w\hbox{-}\lim_{\eps\to 0}}
\mathcal{L}(\frac{-t}{\eps}) \left(
\begin{array}{c}
0 \\  \frac{1}{2}\frac{1}{\sqrt{\rho_{0}}}
\mathbb{H}^{\perp}_{\rho_{0}}[B_{2}(V_{1},V_{2})+B_{2}(V_{2},V_{1})]
\end{array}
\right)=\mathcal{Q}_{2}(V_{1},V_{2})\,.
\end{array}\end{equation}
holds for $V_{1}\in L^{q}([0,T];H^{s}(\T^{n}))$ and $V_{2}\in L^{p}([0,T];H^{-s}(\T^{n}))$.
\end{proposition}

It is also possible to extend the above proposition to the case when  $V_{2}$ in the left hand side is replaced by a sequence $V_{2}^{\eps}$ such that $V_{2}^{\eps}$ converges to $V_{2}$ in $L^{p}([0,T];H^{-s}(\T^{n}))$.

\begin{proposition}\label{pr2}
For any $u$ satisfying $\nabla\cdot(\rho_{0}u)=0$, and for all vectors $V$ and $V_{j}$, $j=1,2$, we have
\begin{equation}\label{pe.2.c}
\int_{\T^{n}}\mathcal{Q}_{1}(u,V)\cdot Vdx=0\,,
\end{equation}
\begin{equation}\label{pe.2.d}
\int_{\T^{n}}\mathcal{Q}_{1}(u,V_{1})\cdot V_{2}+\mathcal{Q}_{1}(u,V_{2})\cdot V_{1}dx=0\,.
\end{equation}
\end{proposition}
\noindent{\it Proof.} From the eigenfunction expansion (\ref{3.8})  of $V$ and spectral decomposition (\ref{3.9}) of $\mathcal{L}(\frac{t}{\eps})V$ we can represent $B_{1}$ defined by (\ref{3.14}) as
\begin{align*}
B_{1}(u,V)&=\displaystyle\sum_{j=1}^{\infty}a_{j}^{+}\frac{i}{\sqrt{\ka_{j}}}e^{i\sqrt{\ka_{j}}t/\eps}
\hbox{div}\big(\rho_{0}u\otimes\nabla\chi_{j}
+\nabla\chi_{j}\otimes\rho_{0}u\big)\\
&\phantom{xx}{}-\sum_{j=1}^{\infty}a_{j}^{-}\frac{i}{\sqrt{\ka_{j}}}
e^{-i\sqrt{\ka_{j}}t/\eps}\hbox{div}\big(\rho_{0}u\otimes\nabla\chi_{j}
+\nabla\chi_{j}\otimes\rho_{0}u\big)\,.
\end{align*}
Direct computation yields
\begin{align*}
&\phantom{xx}{}\int_{\T^{n}} B_{1}(u,V)\cdot
\mathcal{L}_{2}\Big(\frac{t}{\eps}\Big)V \frac{dx}{\sqrt{\rho_{0}}}\\
&=\sum_{j=1}^{\infty}\frac{2a_{j}^{+}a_{j}^{-}}{\ka_{j}}
{\displaystyle\int_{\T^{n}}}\hbox{div}\big(\rho_{0}u\otimes\nabla\chi_{j}
+\nabla\chi_{j}\otimes\rho_{0}u\big)\cdot\nabla\chi_{j}\, dx + R_\eps\,.
\end{align*}
where $R_\eps$ is the oscillating term.
Taking the limit $\eps\to 0$ we have obtained the relation
\begin{equation}\begin{array}{l}
\displaystyle\int_{\T^{n}}\mathcal{Q}_{1}(u,V)\cdot Vdx={\displaystyle\lim_{\eps\to 0}\int_{\T^{n}}}B_{1}(u,V)\cdot\mathcal{L}_{2}
\Big(\frac{t}{\eps}\Big)V\frac{dx}{\sqrt{\rho_{0}}}=0\,.
\end{array}\end{equation}
Then by symmetry, we can conclude
$$
{\displaystyle\int_{\T^{n}}}\mathcal{Q}_{1}(u,V_{1})\cdot V_{2}dx
+{\displaystyle\int_{\T^{n}}}\mathcal{Q}_{1}(u,V_{2})\cdot V_{1}dx=0\,.
$$
\qed

\begin{proposition}\label{pr3}
For any $u$ satisfying $\nabla\cdot(\rho_{0}u)=0$, and for all vectors $V$ and $V_{j}$, $j=1,2$, we have
\begin{equation}
\lim_{\eps\to 0}{\displaystyle\int_{0}^{t}\int_{\T^{n}}}B_{2}(V_{1},V_{2})
\cdot u dxds=0\,,
\end{equation}
\begin{equation}\begin{array}{l}\label{osc1}
{\displaystyle\int_{\T^{n}}}\mathcal{Q}_{2}(V,V)\cdot V dx=0\,,
\end{array}\end{equation}
and
\begin{equation}\begin{array}{l}\label{osc2}
{\displaystyle\int_{\T^{n}}}\mathcal{Q}_{2}(V_{1},V_{1})\cdot V_{2}+2\mathcal{Q}_{2}(V_{1},V_{2})\cdot V_{1}dx=0\,.
\end{array}\end{equation}
\end{proposition}
\noindent{\it Proof.} Straightforward computation gives
\begin{align*}
&\phantom{xx}{}\hbox{div}\Big(\mathcal{L}_{2}\Big(\frac{t}{\eps}\Big)V\otimes
\mathcal{L}_{2}\Big(\frac{t}{\eps}\Big)V\Big)\\
&=\frac{1}{2}\rho_{0}\nabla
\Big(\Big|\frac{\mathcal{L}_{2}(\frac{t}{\eps})V}{\sqrt{\rho_{0}}}\Big|^{2}\Big)
+\frac{1}{\sqrt{\rho_{0}}}\hbox{div}\Big(\sqrt{\rho_{0}}\mathcal{L}_{2}\Big(\frac{t}{\eps}\Big)V\Big)\mathcal{L}_{2}\Big(\frac{t}{\eps}\Big)V\\
&=\frac{1}{2}\rho_{0}\nabla
\Big(\Big|\frac{\mathcal{L}_{2}(\frac{t}{\eps})V}{\sqrt{\rho_{0}}}\Big|^{2}\Big)
-\eps\frac{1}{\sqrt{\rho_{0}}}\pa_{t}\Big(\mathcal{L}_{1}\Big(\frac{t}{\eps}\Big)V\Big)\mathcal{L}_{2}\Big(\frac{t}{\eps}\Big)V\\
&=\frac{1}{2}\rho_{0}\nabla
\Big(\Big|\frac{\mathcal{L}_{2}(\frac{t}{\eps})V}{\sqrt{\rho_{0}}}\Big|^{2}\Big)
-\eps\frac{1}{\sqrt{\rho_{0}}}\pa_{t}\Big(\mathcal{L}_{1}\Big(\frac{t}{\eps}\Big)V\mathcal{L}_{2}\Big(\frac{t}{\eps}\Big)V\Big)
+\eps\frac{1}{\sqrt{\rho_{0}}}\mathcal{L}_{1}\Big(\frac{t}{\eps}\Big)
  V\pa_{t}\Big(\mathcal{L}_{2}\Big(\frac{t}{\eps}\Big)V\Big)\\
  &=\frac{1}{2}\rho_{0}\nabla
\Big(\Big|\frac{\mathcal{L}_{2}(\frac{t}{\eps})V}{\sqrt{\rho_{0}}}\Big|^{2}\Big)
-\frac{1}{2}\nabla\Big(\Big|\mathcal{L}_{1}\Big(\frac{t}{\eps}\Big)V\Big|^{2}\Big)-\eps\frac{1}{\sqrt{\rho_{0}}}\pa_{t}
\Big(\mathcal{L}_{1}\Big(\frac{t}{\eps}\Big)V\mathcal{L}_{2}\Big(\frac{t}{\eps}\Big)V\Big)\,.
\end{align*}
Therefore the bilinear form $B_{2}(V,V)$ defined by (\ref{3.15}) becomes
\begin{equation*}\begin{array}{l}
B_{2}(V,V)=\frac{1}{2}\rho_{0}\nabla
\big(\big|\frac{\mathcal{L}_{2}(\frac{t}{\eps})V}{\sqrt{\rho_{0}}}\big|^{2}\big)
-\eps\frac{1}{\sqrt{\rho_{0}}}\pa_{t}(\mathcal{L}_{1}(\frac{t}{\eps})V\mathcal{L}_{2}(\frac{t}{\eps})V)\,,
\end{array}\end{equation*}
and hence
\begin{align*}
{\displaystyle\int_{0}^{t}\int_{\T^{n}}}B_{2}(V,V)\cdot udxds
&=-\eps{\displaystyle\int_{\T^{n}}}\frac{1}{\sqrt{\rho_{0}}}\Big(\mathcal{L}_{1}\Big(\frac{s}{\eps}\Big)V\mathcal{L}_{2}\Big(\frac{s}{\eps}\Big)V\Big)\cdot
udx\Big|_{0}^{t}\\
&\phantom{xx}{}+\eps{\displaystyle\int_{0}^{t}\int_{\T^{n}}}\frac{1}{\sqrt{\rho_{0}}}
\Big(\mathcal{L}_{1}\Big(\frac{s}{\eps}\Big)V\mathcal{L}_{2}\Big(\frac{s}{\eps}\Big)V\Big)\cdot
\pa_{t}udxds
\end{align*}
which converges to 0 as $\eps$ tends to 0.
Thus by symmetry, we have
\begin{equation*}
\lim_{\eps\to 0}{\displaystyle\int_{0}^{t}\int_{\T^{n}}}B_{2}(V_{1},V_{2})
\cdot u dxds= 0\,.
\end{equation*}
Next, for (\ref{osc1}), we need to study the behavior of $B_{2}(V,V)\cdot
\mathcal{L}_{2}(\frac{t}{\eps})V$.  Let us define the resonance triple $(\ga_{j},\ga_{l},\ga_{m})$ by $\ga_{j}+\ga_{l}+\ga_{m}=0$, where
$$
\ga_{s}=\sqrt{\ka_{s}}\quad \hbox{or}\quad  -\sqrt{\ka_{s}}, \quad s=j, l, m\,.
$$
If the eigenvalues of $L$ remain away from the resonant set then it is easy to obtain (\ref{osc1}) by a stationary phase argument.
Therefore, we only need to focus on the case when
the eigenvalues of $L$ have resonance triples. From (\ref{3.15}) we can rewrite $B_{2}(V,V)$ as
\begin{align*}
B_{2}(V,V)=\frac{1}{2}\rho_{0}\nabla
\Big(\Big|\frac{\mathcal{L}_{2}(\frac{t}{\eps})V}{\sqrt{\rho_{0}}}
     \Big|^{2}\Big)
+\frac{1}{\sqrt{\rho_{0}}}\hbox{div}\Big(\sqrt{\rho_{0}}\mathcal{L}_{2}\Big(\frac{t}{\eps}\Big)V\Big)\mathcal{L}_{2}\Big(\frac{t}{\eps}\Big)V+
\frac{1}{2}\nabla\Big(\Big|\mathcal{L}_{1}\Big(\frac{t}{\eps}\Big)V\Big|^{2}\Big)\,,
\end{align*}
and we will calculate the non-oscillating terms corresponding to the resonance triple $(\ga_{j},\ga_{l},\ga_{m})$.
To proceed, we define
$$
C_{jlm}=a_{j}^{{\rm
sgn}(\ga_{j})}a_{l}^{{\rm sgn}(\ga_{l})}a_{m}^{{\rm sgn}(\ga_{m})}\,,
$$
where ($k=j,l,m$)
\begin{equation*}
a_{k}^{{\rm sgn}(\ga_{k})}=
\left\{\begin{array}{l}
 a_k^+\quad \hbox{when}\quad \ga_{k}>0 \,,
\\ \\
a_k^-\quad \hbox{when}\quad \ga_{k}<0 \,.
\\
\end{array}
\right.\end{equation*}
For any vector $V$ given by eigenfunction expansion (\ref{3.8}), we have
\begin{equation}
\displaystyle\int_{\T^{n}}B_{2}(V,V)\cdot
\mathcal{L}_{2}\Big(\frac{t}{\eps}\Big)V\frac{dx}{\sqrt{\rho_{0}}}
=(I_{1}+I_{2}+I_{3})(\chi_{j},\chi_{l},\chi_{m})+\hbox{other terms}\,,
\end{equation}
where
\begin{align*}
I_{1}&=\frac{-i (C_{jlm}-\overline{C_{jlm}})}{\ga_{j}\ga_{l}\ga_{m}} \int_{\T^{n}}
\rho_{0}\nabla\big(\nabla\chi_{l}\cdot\nabla\chi_{m}\big)\nabla\chi_{j}\\
&\phantom{xx}{}+\rho_{0}\nabla\big(\nabla\chi_{j}\cdot\nabla\chi_{m}\big)\nabla\chi_{l}+\rho_{0}\nabla\big(\nabla\chi_{j}\cdot\nabla\chi_{l}\big)\nabla\chi_{m}
dx\,,
\end{align*}
\begin{align*}
 I_{2}&=\frac{-2i (C_{jlm}-\overline{C_{jlm}})}{\ga_{j}\ga_{l}\ga_{m}} \int_{\T^{n}}
\big(\nabla\chi_{l}\cdot\nabla\chi_{m}\big)\nabla\cdot(\rho_{0}\nabla\chi_{j})\\
&\phantom{xx}{}+\big(\nabla\chi_{j}\cdot\nabla\chi_{m}\big)\nabla\cdot(\rho_{0}\nabla\chi_{l})
+\big(\nabla\chi_{j}\cdot\nabla\chi_{l}\big)\nabla\cdot(\rho_{0}\nabla\chi_{m})
dx\,,
\end{align*}
thus
\begin{equation*}
\displaystyle I_{1}+I_{2}=\frac{i (C_{jlm}-\overline{C_{jlm}})}{\ga_{j}\ga_{l}\ga_{m}} \int_{\T^{n}}
\ga^{2}_{j}\chi_{j}\nabla\chi_{l}\cdot\nabla\chi_{m}
+\ga^{2}_{l}\chi_{l}\nabla\chi_{j}\cdot\nabla\chi_{m}+\ga^{2}_{m}\chi_{m}\nabla\chi_{j}\cdot\nabla\chi_{l}dx
\,,
\end{equation*}
and
\begin{align*}
I_{3}&=i{\displaystyle (C_{jlm}-\overline{C_{jlm}})\int_{\T^{n}}}\frac{1}{\ga_{m}}\nabla(\chi_{j}\chi_{l})\cdot\nabla\chi_{m}
+\frac{1}{\ga_{j}}\nabla(\chi_{l}\chi_{m})\cdot\nabla\chi_{j}
+\frac{1}{\ga_{l}}\nabla(\chi_{j}\chi_{m})\cdot\nabla\chi_{l}dx\,.
\end{align*}
Using the resonance condition $\ga_{j}+\ga_{l}+\ga_{m}=0$, we have
\begin{equation}\begin{array}{l}
I_{1}+I_{2}+I_{3}
=0\,,
\end{array}\end{equation}
and hence
\begin{equation}\begin{array}{l}
\displaystyle\int_{\T^{n}}\mathcal{Q}_{2}(V,V)\cdot Vdx={\displaystyle\lim_{\eps\to 0}\int_{\T^{n}}}\frac{1}{\sqrt{\rho_{0}}}B_{2}(V,V)
\cdot\mathcal{L}_{2}\Big(\frac{t}{\eps}\Big)V dx=0\,.
\end{array}\end{equation}
This proved (\ref{osc1}). Furthermore, we can apply (\ref{osc1}) to $V_{1}+ x V_{2}$, then
\begin{equation}\begin{array}{l}
{\displaystyle\int_{\T^{n}}}\mathcal{Q}_{2}(V_{1}+ x V_{2},V_{1}+ x V_{2})\cdot (V_{1}+ x V_{2}) dx=0\,,
\end{array}\end{equation}
and (\ref{osc2}) follows by equating the terms of degree 1.
\qed

The following property is simple, but useful for our discussion.
\begin{proposition}\label{pro7}
For any $V$ and weighted divergence free $u_{j}$, $\nabla\cdot(\rho_{0}u_{j})=0$, $j=1,2$, as $\eps\to 0$, we have
\begin{equation}
\int_{\T^{n}}B_{1}(u_{1},V)\cdot u_{2}dx\to 0\,,
\end{equation}
\begin{equation}\begin{array}{l}
\displaystyle\int_{\T^{n}}\frac{1}{\sqrt{\rho_{0}}}\,{\rm
div}(\rho_{0}u_{1}\otimes u_{2})\cdot
\mathcal{L}_{2}\Big(\frac{t}{\eps}\Big)Vdx\to 0\,.
\end{array}\end{equation}
\end{proposition}

\noindent{\bf Step 3.} {\it The modulated energy functional and uniform estimates.\it} Let $V^{0}$ be the solution of the oscillating part:
\begin{equation}\label{pe.3.a}
\pa_{t}V^{0}+\mathcal{Q}_{1}(v,V^{0})+\mathcal{Q}_{2}(V^{0},V^{0})=0
\end{equation}
with initial condition
\begin{equation}\label{pe.3.b}
V^{0}(x,0)=V^0_{\rm in}(x)=(\varphi_{0},\sqrt{\rho_{0}}w_{0})^{t}\,,
\end{equation}
where $v$ is the strong solution of the anelastic system (\ref{pe.0.a}).

The local existence of classical solutions of the oscillating system $(\ref{pe.3.a})$--$(\ref{pe.3.b})$ follows the line of \cite{[Mas01a]} with modification. However, for completeness we will give a proof by the standard energy estimate.

\begin{proposition}\label{pro77}
Let $s>\frac{n}{2}+3$ and $v\in C([0,T]; H^{s}(\T^{n}))$ be the local classical solution of the anelastic system $(\ref{pe.0.a})$. Given $V^0_{\rm in}(x)\in H^{s}(\T^{n})$ there exists $T^{*}>0$ and a function $V^{0}\in C([0,T^{*}]; H^{s}(\T^{n}))$ that solves the initial value problem of the oscillating system $(\ref{pe.3.a})$--$(\ref{pe.3.b})$. Moreover, it satisfies the energy equality:
\begin{equation}\begin{array}{c}
{\displaystyle\frac{1}{2}\int_{\T^{n}}}|V^{0}|^{2}dx={\displaystyle\frac{1}{2}\int_{\T^{n}}}|\mathcal{L}(\frac{t}{\eps})V^{0}|^{2}dx
={\displaystyle\frac{1}{2}\int_{\T^{n}}}|V^{0}_{\rm in}|^{2}dx\,,
\end{array}\end{equation}
and hence
\begin{equation}\begin{array}{c}
{\displaystyle\frac{d}{dt}\int_{\T^{n}}}e dx=0\,,
\end{array}\end{equation}
where the energy density defined by
\begin{equation}\begin{array}{l}
e={\displaystyle\frac{1}{2}}\rho_{0}|v|^{2}
+{\displaystyle\frac{1}{2}}|\mathcal{L}(\frac{t}{\eps})V^{0}|^{2}\,.
\end{array}\end{equation}
\end{proposition}
\noindent{\it Proof.}
We are just going to give the priori bounds we can derive for this system. The existence result is then easily deduced by solving some approximated systems. For instance, we can project the system to
\begin{equation}
\mathcal{P}_{M}=\hbox{span}\bigg\{\pm i\sqrt{\ka_{j}}\,,\left(
\begin{array}{c}
\chi_{j} \\ \frac{\pm i}{\sqrt{\ka_{j}}}\sqrt{\rho_{0}}\nabla\chi_{j}
\end{array}
\right)\bigg\}\,,\quad \hbox{where} \quad \ka_{j}<M\,.
\end{equation}
Then, we have just to take the limit $M\to \infty$ and use a compactness method to pass to the limit.

Now, we turn to the proof of the a priori bounds. By Proposition \ref{pr2}, we have
$$
\int_{\T^{n}}\mathcal{Q}_{1}(v,V^{0})\cdot V^{0}dx=0
\quad\hbox{and}\quad \int_{\T^{n}}\mathcal{Q}_{2}(V^{0},V^{0})\cdot V^{0}dx=0\,,
$$
then the standard energy estimate gives the $L^2$-estimate
$$
\frac{d}{dt}\|V^{0}\|^{2}_{L^{2}(\T^{n})}=0\,.
$$
The higher energy estimate is obtained in the similar way. We will discuss $\mathcal{Q}_{1}$ first. Indeed, for any $j\leq s$, we have
\begin{equation}\label{loc.1a}
\displaystyle\int_{\T^{n}}\pa^{j}_{x}\mathcal{Q}_{1}(v,V^{0})\cdot \pa^{j}_{x}V^{0}dx
=
{\displaystyle\lim_{\eps\to 0}}\int_{\T^{n}}\mathcal{L}\Big(\frac{-t}{\eps}\Big)
\pa^{j}_{x}\left(
\begin{array}{c}
0 \\  \frac{1}{\sqrt{\rho_{0}}}\mathbb{H}^{\perp}_{\rho_{0}}B_{1}(v,V^{0})
\end{array}
\right)\cdot \pa^{j}_{x}V^{0}dx\,.
\end{equation}
Thus we only need to estimate
\begin{equation}\begin{array}{l}\label{loc.1}
\displaystyle{\lim_{\eps\to 0}\int_{\T^{n}}}B_{1}(v,\pa^{j}_{x}V^{0})\cdot\mathcal{L}_{2}
\Big(\frac{t}{\eps}\Big)\pa^{j}_{x}V^{0}\frac{dx}{\sqrt{\rho_{0}}}
\end{array}\end{equation}
since all the other terms of (\ref{loc.1a}) can be controlled by low order energy norm. However, it is easy to see that (\ref{loc.1}) is zero by Proposition \ref{pr2} and we have
\begin{equation}
\Big|\int_{\T^{n}}\pa^{j}_{x}\mathcal{Q}_{1}(v,V^{0})\cdot \pa^{j}_{x}V^{0}dx\Big|\leq C_{1}\|V^{0}\|^{2}_{H^{j}(\T^{n})}\,.
\end{equation}
For $\mathcal{Q}_{2}$, note that
\begin{equation}\label{qq1}
\displaystyle\int_{\T^{n}}\pa^{j}_{x}\mathcal{Q}_{2}(V^{0},V^{0})\cdot \pa^{j}_{x}V^{0}dx
=
{\displaystyle\lim_{\eps\to 0}}\int_{\T^{n}}\mathcal{L}\Big(\frac{-t}{\eps}\Big)
\pa^{j}_{x}\left(
\begin{array}{c}
0 \\  \frac{1}{\sqrt{\rho_{0}}}\mathbb{H}^{\perp}_{\rho_{0}}B_{2}(V^{0},V^{0})
\end{array}
\right)\cdot \pa^{j}_{x}V^{0}dx\,.
\end{equation}
Similarly, by Proposition \ref{pr3} and Sobolev inequality,  we also have the following estimate
\begin{align*}
&\phantom{xx}{}\displaystyle{\lim_{\eps\to 0}\int_{\T^{n}}}B_{2}(V^{0},\pa^{j}_{x}V^{0})\cdot\mathcal{L}_{2}
\Big(\frac{t}{\eps}\Big)\pa^{j}_{x}V^{0}\frac{dx}{\sqrt{\rho_{0}}}\\
&=-\displaystyle{\lim_{\eps\to 0}\frac{1}{2}\int_{\T^{n}}}B_{2}(\pa^{j}_{x}V^{0},\pa^{j}_{x}V^{0})\cdot\mathcal{L}_{2}
\Big(\frac{t}{\eps}\Big)V^{0}\frac{dx}{\sqrt{\rho_{0}}}\\
&\leq C(\rho_{0})\|V^{0}\|^{2}_{H^{j}(\T^{n})}\|\nabla_{x} V^{0}\|_{L^{\infty}(\T^{n})}\leq C(\rho_{0})\|V^{0}\|^{2}_{H^{j}(\T^{n})}\|V^{0}\|_{H^{s}(\T^{n})}\,.
\end{align*}
Moreover, the other terms of (\ref{qq1}) can be bounded by
$$
\sum_{s+r=j+1, j\geq r, s\geq 1}\|\pa^{r}_{x}V^{0}\|_{L^{4}(\T^{n})}\|\pa^{s}_{x}V^{0}\|_{L^{4}(\T^{n})}\|\pa^{j}_{x}V^{0}\|_{L^{2}(\T^{n})}\,.
$$
Then, using Gagliardo-Nirenberg inequality, we deduce
$$
\|\pa^{r}_{x}V^{0}\|_{L^{4}(\T^{n})}\leq C \|\pa_{x}V^{0}\|^{1-\theta}_{L^{\infty}(\T^{n})} \|\pa^{j-1}_{x}V^{0}\|^{\theta}_{L^{2}(\T^{n})}\,,
$$
where $\frac{1}{4}-(r-1)=\theta(\frac{1}{2}-(j-1))$. And, since $r+s=j+1$, we also have
$$
\|\pa^{s}_{x}V^{0}\|_{L^{4}(\T^{n})}\leq C \|\pa_{x}V^{0}\|^{\theta}_{L^{\infty}(\T^{n})} \|\pa^{j-1}_{x}V^{0}\|^{1-\theta}_{L^{2}(\T^{n})}\,.
$$
Combining all the above energy estimates we have
$$
\frac{d}{dt}\|V^{0}\|^{2}_{H^{s}(\T^{n})}\leq C_{1}\|V^{0}\|^{2}_{H^{s}(\T^{n})}+C_{2}\|V^{0}\|^{3}_{H^{s}(\T^{n})}\,.
$$
Once the energy estimate is obtained, the local existence of classical solutions follows immediately by the standard method and we omit the detail.
\qed

We can define the modulated energy $H^{\eps}(t)$ of (\ref{int.1.a}) as
\begin{equation}
 H^{\eps}(t)= {\displaystyle\frac{1}{2}\int_{\T^{n}}}
\Big|\Big(\eps^{\alpha}\nabla-i\big[v+\frac{1}{\sqrt{\rho_{0}}}\mathcal{L}_{2}\Big(\frac{t}{\eps}\Big)V^{0}\big]\Big)\psi^{\eps}\Big|^{2}dx
+ {\displaystyle\frac{1}{2}\int_{\T^{n}}}\Big|\varphi^{\eps}
-\mathcal{L}_{1}\Big(\frac{t}{\eps}\Big)V^{0}\Big|^{2}dx
\end{equation}
where $\mathcal{L}_1$ and $\mathcal{L}_2$ are given by (\ref{3.10}).
For simplicity, we define
\begin{equation}
\varpi=v+\frac{1}{\sqrt{\rho_{0}}}\mathcal{L}_{2}\Big(\frac{t}{\eps}\Big)V^{0}
\end{equation}
and rewrite the modulated energy $H^{\eps}(t)$ as
\begin{align*}
H^{\eps}(t)&={\displaystyle\int_{\T^{n}}}(e^{\eps}
+e) dx+{\displaystyle\int_{\T^{n}}\frac{1}{2}}
(\rho^{\eps}-\rho_{0})|\varpi|^{2}dx\\
&\phantom{xx}{}-{\displaystyle\int_{\T^{n}}}J^{\eps}\cdot vdx
-{\displaystyle\int_{\T^{n}}}\sqrt{\rho_{0}}\nabla
w^{\eps}\cdot\mathcal{L}_{2}\Big(\frac{t}{\eps}\Big)V^{0}dx
-{\displaystyle\int_{\T^{n}}}\mathcal{L}_{1}\Big(\frac{t}{\eps}\Big)V^{0}\cdot\varphi^{\eps}dx\,.
\end{align*}
Using the identity
\begin{equation}\begin{array}{c}
 \pa_{t}(\mathcal{L}(\frac{t}{\eps})V^{0})=
 \mathcal{L}(\frac{t}{\eps})\pa_{t}V^{0}+\frac{1}{\eps}\mathcal{L}(\frac{t}{\eps})LV^{0}
\,,
\end{array}\end{equation}
we have the evolution of the modulated energy
\begin{align}\label{3.47}
{\displaystyle\frac{d}{dt}H^{\eps}}(t)&=
{\displaystyle\frac{d}{dt}\int_{\T^{n}}\frac{1}{2}}(\rho^{\eps}-\rho_{0})
|\varpi|^{2}dx-{\displaystyle\int_{\T^{n}}}J^{\eps}\cdot \pa_{t}vdx\\
&\phantom{xx}{}-{\displaystyle\int_{\T^{n}}}V^{\eps}\cdot\pa_{t}V^{0}dx+A_1
+A_2+A_3\nonumber
\end{align}
where $A_1$, $A_2$ and $A_3$ are given respectively by
\begin{equation}\begin{array}{l}
 A_{1}=-{\displaystyle\frac{\eps^{2\alpha}}{2} \int_{\T^{n}}}
 \big(\nabla\psi^{\eps}\otimes\nabla\overline{\psi^\eps}
 +\nabla\overline{\psi^\eps}\otimes\nabla\psi^{\eps}
\big):\nabla  \varpi dx\,,
\\[2mm]
 A_{2}=-{\displaystyle\int_{\T^{n}}\frac{1}{2}}(\varphi^{\eps})^{2}\,
\hbox{div}  \varpi dx
\quad\hbox{and}\quad
 A_{3}={\displaystyle\frac{\eps^{2\alpha}}{4}\int_{\T^{n}}}\rho^{\eps}\Delta\,
 \hbox{div} \varpi dx\,.
\end{array}
\end{equation}
We can further rewrite $A_{1}$ as
\begin{align*}
A_{1}&=-{\displaystyle\frac{1}{2}
\int_{\Om}}\Big(\big(\eps^{\alpha}\nabla -i \varpi \big)\psi^{\eps}\otimes\overline{\big(\eps^{\alpha}\nabla -i \varpi \big)\psi^{\eps}}\\
&\phantom{xx}{}+\overline{\big(\eps^{\alpha}\nabla -i \varpi
\big)\psi^{\eps}}\otimes\big(\eps^{\alpha}\nabla -i \varpi
\big)\psi^{\eps}\Big) :\nabla  \varpi
dx+\displaystyle\sum_{i=1}^{6}K_{i}\,,
\end{align*}
with $K_i, i=1,\cdots, 6,$ given by
\begin{align*}
K_{1}&={\displaystyle-\int_{\T^{n}}}\big(v\otimes J^{\eps}+J^{\eps}\otimes v-\rho_{0}v\otimes v\big):\nabla vdx\,,
\\
K_{2}&={\displaystyle\int_{\T^{n}}}\mathcal{L}_{2}\Big(\frac{t}{\eps}\Big)V^{0}\otimes\mathcal{L}_{2}\Big(\frac{t}{\eps}\Big)V^{0}
:\nabla
 \varpi dx\,,
\\
K_{3}&=-{\displaystyle\int_{\T^{n}}}\Big(\frac{1}{\sqrt{\rho_{0}}}\mathcal{L}_{2}\Big(\frac{t}{\eps}\Big)V^{0}\otimes J^{\eps}
 +J^{\eps}\otimes \frac{1}{\sqrt{\rho_{0}}}\mathcal{L}_{2}\Big(\frac{t}{\eps}\Big)V^{0}\Big)
:\nabla
 \varpi dx\,,
\\
 K_{4}&={\displaystyle\int_{\T^{n}}}\Big(\sqrt{\rho_{0}}\mathcal{L}_{2}\Big(\frac{t}{\eps}\Big)V^{0}\otimes v
 +\sqrt{\rho_{0}}v\otimes \mathcal{L}_{2}\Big(\frac{t}{\eps}\Big)V^{0}\Big)
:\nabla
 \varpi dx\,,
\\
K_{5}&={\displaystyle-\int_{\T^{n}}}\Big(v\otimes J^{\eps}+J^{\eps}\otimes v-\rho_{0}v\otimes v\Big):\nabla \Big(\frac{1}{\sqrt{\rho_{0}}}\mathcal{L}_{2}\Big(\frac{t}{\eps}\Big)V^{0}\Big)dx\,.
\end{align*}
and
\begin{align*}
K_{6}&={\displaystyle\int_{\T^{n}}}(\rho^{\eps}-\rho_{0})\big[
\sqrt{\rho_{0}}\mathcal{L}_{2}\Big(\frac{t}{\eps}\Big)V^{0}\otimes v
 +\sqrt{\rho_{0}}v\otimes \mathcal{L}_{2}\Big(\frac{t}{\eps}\Big)V^{0}\big]
:\nabla \displaystyle \varpi \frac{dx}{\rho_{0}}\\
&\phantom{xx}{} +{\displaystyle\int_{\T^{n}}}(\rho^{\eps}-\rho_{0})
\big[\mathcal{L}_{2}\Big(\frac{t}{\eps}\Big)V^{0}\otimes\mathcal{L}_{2}
 (\frac{t}{\eps})V^{0}+\rho_{0}v\otimes v\big]
 :\nabla \displaystyle \varpi \frac{dx}{\rho_{0}}\,.
\end{align*}
Note that $A_{2}$ can be rewritten as
\begin{equation}\begin{array}{l}
 A_{2}=-{\displaystyle\int_{\T^{n}}\frac{1}{2}}|\varphi^{\eps}-
 \mathcal{L}_{1}(\frac{t}{\eps})V^{0}|^{2}
\hbox{div}  \varpi dx+K_{7}+K_{8}\,,
\end{array}\end{equation}
where $K_7$ and $K_8$ are given by
\begin{align*}
 K_{7}&={\displaystyle\int_{\T^{n}}}{\displaystyle\frac{1}{2}}\Big|
 \mathcal{L}_{1}\Big(\frac{t}{\eps}\Big)V^{0}\Big|^{2}
\hbox{div}  \varpi dx\,,\\
 K_{8}&={\displaystyle -\int_{\T^{n}}}\varphi^{\eps}\cdot
 \mathcal{L}_{1}\Big(\frac{t}{\eps}\Big)V^{0}
\hbox{div}  \varpi dx\,.
\end{align*}
Using the following identities
\begin{align*}
-\int_{\T^{n}}\big({v}\otimes {J}^{\eps}\big):\nabla {v}dx
&=\int_{\T^{n}}\frac{1}{2}|{v}|^{2}\nabla\cdot
J^{\eps}dx\,,\\
-\int_{\T^{n}}(J^{\eps}\otimes {v}\big):\nabla {v}dx
&=-\int_{\T^{n}}\Big[\big({v}\cdot\nabla\big){v}\Big]\cdot J^{\eps}dx\,,
\end{align*}
and
\begin{align*}
 \int_{\T^{n}}\big(\rho_{0}{v}\otimes {v}\big):\nabla {v}dx
&=\int_{\T^{n}}\big[\big(\rho_{0}{v}\cdot\nabla\big){v}\big]\cdot v dx\\
\displaystyle&=\frac{1}{2}\int_{\T^{n}}\big(\rho_{0}{v}\cdot\nabla\big)|v|^{2} dx=0\,,
\end{align*}
we can rewrite $K_{1}$ as
\begin{equation}
\displaystyle K_{1}=\int_{\T^{n}}\frac{1}{2}|v|^{2}\nabla\cdot
J^{\eps}
-\big[\big(v\cdot\nabla\big)v\big]\cdot J^{\eps}dx\,.
\end{equation}
For $K_{2}+K_{7}$, by Proposition \ref{pro7}, we have
\begin{align*}
K_{2}+K_{7}
 &=-\displaystyle\int_{\T^{n}}B_{2}(V^{0},V^{0})
\cdot
( \rho_{0}\varpi )\frac{dx}{\rho_{0}}\\
&=-{\displaystyle\int_{\T^{n}}}B_{2}(V^{0},V^{0})
\cdot vdx+o_{\eps}(1)\,.
\end{align*}
For $K_{3}+K_{8}$, by proposition \ref{pr2} and \ref{pro7}, we have
\begin{align*}
&\phantom{xx}{}K_{3}+K_{8}\\
 &={\displaystyle\int_{\T^{n}}}
 \hbox{div}\Big[\frac{1}{\sqrt{\rho_{0}}}\mathcal{L}_{2}\Big(\frac{t}{\eps}\Big)V^{0}\otimes J^{\eps}
 +J^{\eps}\otimes \frac{1}{\sqrt{\rho_{0}}}\mathcal{L}_{2}\Big(\frac{t}{\eps}\Big)V^{0}\Big]
\cdot \displaystyle( \rho_{0}\varpi )\frac{dx}{\rho_{0}}\\
&\phantom{xx}{} +{\displaystyle\int_{\T^{n}}}
 \nabla\big(\varphi^{\eps}\cdot\mathcal{L}_{1}\Big(\frac{t}{\eps}\Big)V^{0}\big)
 \cdot \displaystyle( \rho_{0}\varpi )\frac{dx}{\rho_{0}}\\
& ={\displaystyle\int_{\T^{n}}}
 \hbox{div}\Big[\sqrt{\rho_{0}}\mathcal{L}_{2}\Big(\frac{t}{\eps}\Big)V^{0}\otimes \nabla w^{\eps}
 +\nabla w^{\eps}\otimes \sqrt{\rho_{0}}\mathcal{L}_{2}\Big(\frac{t}{\eps}\Big)V^{0}\Big]
\cdot\displaystyle ( \rho_{0}\varpi )\frac{dx}{\rho_{0}} \\
&\phantom{xx}{}+{\displaystyle\int_{\T^{n}}}
\hbox{div}\Big[\sqrt{\rho_{0}}\mathcal{L}_{2}\Big(\frac{t}{\eps}\Big)V^{0}\otimes u^{\eps}
 +u^{\eps}\otimes \sqrt{\rho_{0}}\mathcal{L}_{2}\Big(\frac{t}{\eps}\Big)V^{0}\Big]
\cdot \displaystyle( \rho_{0}\varpi )\frac{dx}{\rho_{0}}\\
 &\phantom{xx}{}+{\displaystyle\int_{\T^{n}}}
 \nabla\big(\varphi^{\eps}\cdot\mathcal{L}_{1}\Big(\frac{t}{\eps}\Big)V^{0}\big)
 \cdot \displaystyle( \rho_{0}\varpi )\frac{dx}{\rho_{0}}\\
&={\displaystyle\int_{\T^{n}}}2B_{2}(V^{0},V^{\eps}) \cdot (
\rho_{0}\varpi
)\displaystyle\frac{dx}{\rho_{0}}+{\displaystyle\int_{\T^{n}}}
B_{1}(u^{\eps},V^{0})
 \cdot ( \rho_{0}\varpi )\displaystyle\frac{dx}{\rho_{0}}\\
\displaystyle&= {\displaystyle\int_{\T^{n}}}2\mathcal{Q}_{2}(V^{0},\overline{V}) \cdot V^{0}dx+{\displaystyle\int_{\T^{n}}}2B_{2}(V^{0},V^{\eps}) \cdot vdx+o_{\eps}(1)\,.
\end{align*}
Similarly, we have
\begin{equation}\begin{array}{l}
 K_{4}=-{\displaystyle\int_{\T^{n}}} B_{1}(v,V^{0})
 \cdot ( \rho_{0}\varpi )\displaystyle\frac{dx}{\rho_{0}}=o_{\eps}(1)\,.
\end{array}\end{equation}
and
\begin{align*}
 K_{5}&={\displaystyle-\int_{\T^{n}}}\big(v\otimes \rho_{0}u^{\eps}+ \rho_{0}u^{\eps}\otimes v-
 \rho_{0}v\otimes v\big):\nabla\Big( \frac{1}{\sqrt{\rho_{0}}}\mathcal{L}_{2}\Big(\frac{t}{\eps}\Big)V^{0}\Big)dx\\
&\phantom{xx}{}-{\displaystyle\int_{\T^{n}}}\big(v\otimes \rho_{0}\nabla w^{\eps}+ \rho_{0}\nabla w^{\eps}\otimes v\big):\nabla \Big( \frac{1}{\sqrt{\rho_{0}}}\mathcal{L}_{2}\Big(\frac{t}{\eps}\Big)V^{0}\Big)dx\\
&={\displaystyle\int_{\T^{n}}} \mathcal{Q}_{1}(v,\overline{V})\cdot V^{0}dx+o_{\eps}(1)\,.
\end{align*}
Finally, it is easy to see $K_{6}=o_{\eps}(1)$, then using the uniform
bound of energy, one can show that $A_{3}=o_{\eps}(1)$.
From (\ref{3.47}) it is natural to define
\begin{equation}\begin{array}{l} W^{\eps}(t)=-{\displaystyle\int_{\T^{n}}\frac{1}{2}}
(\rho^{\eps}-\rho_{0})|\varpi|^{2}dx\,.
\end{array}\end{equation}
which is the measure of density variation. More precisely, it is designed to control the homogeneity of the fluid. We now combine all above computations into (\ref{3.47}) to obtain
\begin{align}\label{3.63}
 &\phantom{xx}{}\frac{d}{dt}\big(H^{\eps}(t)+W^{\eps}(t)\big)\nonumber\\
 &\leq o_{\eps}(1)+H^{\eps}(t)+{\displaystyle\int_{\T^{n}}}
 \big[2B_{2}(V^{0},V^{\eps})-B_{2}(V^{0},V^{0})\big] \cdot vdx\\
 &\phantom{xx}{}+{\displaystyle\int_{\T^{n}}}\big(\mathcal{Q}_{1}(v,\overline{V})\cdot
V^{0}+2\mathcal{Q}_{2}(V^{0},\overline{V}) \cdot V^{0}-V^{\eps}\cdot\pa_{t}V^{0}\big)dx\nonumber\\
&\phantom{xx}{}+\int_{\T^{n}}\Big(\frac{1}{2}|v|^{2}\nabla\cdot
J^{\eps}
+\nabla p\cdot J^{\eps}\Big)dx
\,,\nonumber
\end{align}
Note that by Proposition \ref{pr2} and \ref{pr3}, the second integral on the right hand side of (\ref{3.63}) tends to 0 as $\eps\to 0$;
\begin{align*}
&\phantom{xx}{}{\displaystyle\int_{\T^{n}}}\big(\mathcal{Q}_{1}(v,\overline{V})\cdot
V^{0}+2\mathcal{Q}_{2}(V^{0},\overline{V}) \cdot V^{0}-V^{\eps}\cdot\pa_{t}V^{0}\big)dx\\
&=-{\displaystyle\int_{\T^{n}}}\mathcal{Q}_{1}(v,V^{0})\cdot\overline{V}
+\mathcal{Q}_{2}(V^{0},V^{0}) \cdot \overline{V}+\pa_{t}V^{0}\cdot V^{\eps}dx\to 0
\,,
\end{align*}
and using the equation of continuity the last integral of (\ref{3.63}) can be rewritten as
\begin{align*}
&\phantom{xx}{}\displaystyle\int_{\T^{n}}\Big(\frac{1}{2}|v|^{2}\nabla\cdot
J^{\eps}+\nabla p\cdot J^{\eps}\Big)dx=\displaystyle\int_{\T^{n}}\Big(\frac{1}{2}|v|^{2}-p\Big)\pa_{t}(\rho^{\eps}-\rho_{0})dx\\
&=\displaystyle\frac{d}{dt}\int_{\T^{n}}\Big(\frac{1}{2}|v|^{2}-p\Big)(\rho^{\eps}-\rho_{0})dx
-\int_{\T^{n}}\pa_{t}
\Big(\frac{1}{2}|v|^{2}-p\Big)(\rho^{\eps}-\rho_{0})dx\,.
\end{align*}
Therefore we need to introduce one more correction term $S^{\eps}$ which also measures the density variation defined by
\begin{equation}
\displaystyle
S^{\eps}(t)=-\int_{\T^{n}}\Big(\frac{1}{2}|v|^{2}-p\Big)
(\rho^{\eps}-\rho_{0})dx\,.
\end{equation}
Thus (\ref{3.63}) becomes
\begin{align*}
&\phantom{xx}{}\displaystyle\frac{d}{dt}\big(H^{\eps}(t)+W^{\eps}(t)+S^{\eps}(t)\big)\\
&\leq
o_{\eps}(1)+{\displaystyle\int_{\T^{n}}}\Big[2B_{2}(V^{0},V^{\eps})-B_{2}(V^{0},V^{0})\Big] \cdot vdx+H^{\eps}(t)\,.
\end{align*}
It is easy to show that $W^{\eps}(t)=o_{\eps}(1)$ and
$S^{\eps}(t)=o_{\eps}(1)$, thus
\begin{equation}
H^{\eps}(t)\leq C\big(o_{\eps}(1)+H^{\eps}(0)\big)\,.
\end{equation}
According to Gronwall's inequality we only need to estimate the initial
modulated energy $H^{\eps}(0)$. To proceed, we rewrite the modulated
energy as
\begin{align*}
 H^{\eps}(t)&=
{\displaystyle\frac{\eps^{2\alpha}}{2}
\int_{\T^{n}}}\Big|\nabla\sqrt{\rho_{0}^{\eps}}\,\Big|^{2}dx
+{\displaystyle\frac{1}{2}\int_{\T^{n}}}
\Big|\Big(\eps^{\alpha}\nabla-i\big[v_{0}+\nabla w_{0}\big]\Big)\psi_{0}^{\eps}\Big|^{2}dx\\
&\phantom{xx}{}+ \displaystyle\frac{1}{2}\int_{\T^{n}}
|\varphi_{0}^{\eps}-\varphi_{0}|^{2}dx\,.
\end{align*}
It is sufficient to check the kinetic part which can be rewritten as
\begin{align}\label{3.69}
\displaystyle\int_{\T^{n}}
\Big|\frac{1}{\sqrt{\rho_{0}^{\eps}}}J_{0}^{\eps}
-\sqrt{\rho_{0}^{\eps}}\big(v_{0}+\nabla w_{0}\big)\Big|^{2}dx\,.
\end{align}
By triangle inequality, (\ref{3.69}) can be estimated as
\begin{align}\label{pr.5.e}
&\phantom{xx}{}\Big\|\frac{1}{\sqrt{\rho_{0}^{\eps}}}J_{0}^{\eps}-\sqrt{\rho_{0}^{\eps}}\big(v_{0}+\nabla w_{0}\big)\Big\|_{L^{2}(\T^{n})}\nonumber
\\
\displaystyle&\leq\Big\|\frac{1}{\sqrt{\rho^{\eps}_{0}}}J^{\eps}_{0}-
\frac{1}{\sqrt{\rho_{0}}}J_{0}\Big\|_{L^{2}(\T^{n})}
+\Big\|\big(\sqrt{\rho_{0}}-\sqrt{\rho^{\eps}_{0}}\big)\big(v_{0}+\nabla
w_{0}\big)\Big\|_{L^{2}(\T^{n})}\,.
\end{align}
Note that the first term on the right side of (\ref{pr.5.e}) converges to 0 by assumption (A2). For the second term, using assumption (A1) and an elementary inequality
$$
|\sqrt{x}-\sqrt{a}|^{2}\leq a^{-1}|x-a|^{2}, \qquad x\geq 0,
$$
then employing the fact that $a$ is bounded away from 0, $a\geq c>0$, we have
\begin{align}\label{pr.5.ea}
&\phantom{xx}{}\|(\sqrt{\rho_{0}}-\sqrt{\rho^{\eps}_{0}})\big(v_{0}+\nabla w_{0}\big)\|_{L^{2}(\T^{n})}\nonumber
\\
&\leq
\|v_{0}+\nabla w_{0}\|_{L^{\infty}(\T^{n})}
\|\sqrt{\rho_{0}}-\sqrt{\rho^{\eps}_{0}}\|_{L^{2}(\T^{n})}
\\
&\leq \|v_{0}+\nabla w_{0}\|_{L^{\infty}(\T^{n})}
\|\rho_{0}-\rho^{\eps}_{0}\|_{L^{2}(\T^{n})}\,.\nonumber
\end{align}

Since $\rho_0^\eps\to \rho_0$ in $L^2(\T^n)$, the second term on the  right hand side of (\ref{pr.5.e}) also converges to 0,
this shows $H^{\eps}(0)\to 0$ as $\eps\to 0$. Hence $H^{\eps}(t)\to 0$ as $\eps \to 0$.
It is easy to rewrite the modulated energy as
\begin{align}\label{in.9.a}
H^{\eps}(t)&=
\displaystyle\frac{\eps^{2\alpha}}{2}\int_{\T^{n}}
|\nabla\sqrt{\rho^{\eps}}\,|^{2}dx
+{\displaystyle\frac{1}{2}\int_{\T^{n}}}
\bigg|\varphi^{\eps}-\mathcal{L}_{1}\Big(\frac{t}{\eps}\Big)V^{0}\bigg|^{2}dx\nonumber
\\
&\phantom{xx}{}+{\displaystyle\frac{1}{2}\int_{\T^{n}}
\bigg|\frac{1}{\sqrt{\rho^{\eps}}}\Big(J^{\eps} -\rho^{\eps}
v-\frac{\rho^{\eps}}{\sqrt{\rho_{0}}}
\mathcal{L}_{2}\Big(\frac{t}{\eps}\Big)V^{0}\Big)\bigg|^{2}dx}\,.
\end{align}
From (\ref{in.9.a}) we have
\begin{equation}\begin{array}{lll}\label{in.9.b}
{\displaystyle\frac{1}{2}\int_{\T^{n}}
\bigg|\frac{1}{\sqrt{\rho^{\eps}}}\Big(J^{\eps} -\rho^{\eps}
v-\frac{\rho^{\eps}}{\sqrt{\rho_{0}}}
\mathcal{L}_{2}\Big(\frac{t}{\eps}\Big)V^{0}\Big)\bigg|^{2}dx\to 0\,.}
\end{array}\end{equation}
Therefore we deduce from (\ref{in.9.b}) and H\"older inequality that
\begin{align}\label{in.9.c}
&\phantom{xx}{}{\displaystyle\bigg\|J^{\eps}
-\rho_{0} v-\sqrt{\rho_{0}}\mathcal{L}_{2}
\Big(\frac{t}{\eps}\Big)V^{0}\bigg\|_{L^{4\over3}(\T^{n})}}\nonumber
\\
&{\displaystyle\leq\|\sqrt{\rho^{\eps}}\|_{L^{4}
(\T^{n})}
\bigg\|\frac{1}{\sqrt{\rho^{\eps}}}\Big(J^{\eps} -\rho^{\eps} v-\frac{\rho^{\eps}}{\sqrt{\rho_{0}}}
\mathcal{L}_{2}\Big(\frac{t}{\eps}\Big)V^{0}\Big)\bigg\|_{L^{2}(\T^{n})}}\\
&\phantom{xx}{}{\displaystyle +\|\rho^{\eps}-\rho_{0}\|_{L^{2}
(\T^{n})}\bigg\|v+\frac{1}{\sqrt{\rho_{0}}}
\mathcal{L}_{2}\Big(\frac{t}{\eps}\Big)V^{0}\bigg\|_{L^{4}(\T^{n})}\,.}\nonumber
\end{align}
Thus for any $\varphi(x,t)\in L^{1}([0,T];L^{4}(\T^{n}))$, we have
$$
\lim_{\eps\to 0}\int_{0}^{t}\int_{\T^{n}}\big(J^{\eps}(x,s)-\rho_{0}v(x,s)\big)\cdot\varphi(x,s) dxds=0\,.
$$
i.e.,
$$
J^{\eps}\rightharpoonup \rho_{0}v \quad \hbox{weakly $*$ in} \quad
L^{\infty}\big([0,T];L^{4/3}(\T^{n})\big)\,,
$$
by duality argument. This completes the proof of the main theorem.


\end{document}